\documentclass[12pt]{amsart}
\usepackage{amssymb, graphicx}

\makeatletter 

\newcommand{\ed}{

\end{document}}

\newcommand*\nrotarrowconstructor[2]{%
  \mathrel{\m@th\sbox\z@{$ #1 $}%
    \raisebox{1.3\dp\z@}{%
      \makebox[\wd\z@][c]{%
        \reflectbox{\rotatebox[origin=cB]{90}{$ #2 $}}%
        \kern0.32\wd\z@%
      }}}%
}

\newcommand\nuparrow{\nrotarrowconstructor\uparrow\nrightarrow}
\newcommand{\wins}[2]{{#1\!\uparrow\!{#2}}}
\newcommand{\nwins}[2]{{#1\!\nuparrow\!{#2}}}
\newcommand{\GAME}[1]{{#1}_{\mathsf{G}}}

\newcommand{\I}{\mathbf{I}}
\newcommand{\II}{\mathbf{II}}

\newcommand{\cl}[1]{\overline{#1}}
\newcommand{\Union}{\bigcup}

\newcommand{\bbA}{\mathbb{A}}

\newcommand{\seq}[1]{\{#1\}_{n\in\N}}

\newcommand{\cA}{\mathcal{A}}
\newcommand{\cC}{\mathcal{C}}

\newcommand{\cU}{\mathcal{U}}

\long\def\forget#1\forgotten{}

\newcommand{\Impl}{\Rightarrow}

\newcommand{\N}{{\mathbb N}}

\newcommand{\R}{\mathbb R}

\newcommand{\U}{\bigcup}

\newcommand{\nin}{\notin}
\newcommand{\cat}{\hat{\ }}
\newcommand{\sbst}{\subseteq}
\newcommand{\spst}{\supseteq}
\newcommand{\sm}{\setminus}

\newcommand{\A}{\forall}

\long\def\note#1\endnote%
{\ifNote{\par\medskip\tt Note: #1}\medskip\par
                         \else{}
                         \fi}

\newtheorem{thm}{Theorem}[section]
\newtheorem{cor}[thm]{Corollary}
\newtheorem{conj}[thm]{Conjecture}
\newtheorem{lem}[thm]{Lemma}
\newtheorem{prop}[thm]{Proposition}
\newtheorem{prob}[thm]{Problem}

\theoremstyle{remark}

\newtheorem{rem}[thm]{Remark}

\theoremstyle{definition}
\newtheorem{defn}[thm]{Definition}

\newcommand{\be}{\begin{enumerate}}
\newcommand{\ee}{\end{enumerate}}
\newcommand{\bi}{\begin{itemize}}
\newcommand{\ei}{\end{itemize}}

\newcommand{\bpf}{\begin{proof}}
\newcommand{\epf}{\end{proof}}

\title{Null sets and games in Banach spaces}

\author{Jakub Duda}
\address[Jakub Duda]{Department of Mathematics,
Weizmann Institute of Science, Rehovot 76100, Israel.}
\curraddr{\v{C}EZ, a.s., Duhov\'a 2/1444, 140 53 Praha 4, Czech Republic}
\email{jakub.duda@gmail.com}

\author{Boaz Tsaban}
\address[Boaz Tsaban]{
Department of Mathematics, Bar-Ilan University,
Ramat-Gan 52900, Israel;
and
Department of Mathematics,
Weizmann Institute of Science, Rehovot 76100, Israel.}
\email{tsaban@math.biu.ac.il}


\begin{document}

\begin{abstract}
The notion of Aronszajn-null sets generalizes the notion of
Lebesgue measure zero in the Euclidean space to infinite
dimensional Banach spaces. We present a game-theoretic approach to
Aronszajn-null sets, establish its basic properties, and discuss
some ensuing open problems.
\end{abstract}

\maketitle

\section{Motivation}

Aronszajn null sets were introduced by Aronszajn in the context of
studying almost-everywhere differentiability of Lipschitz mappings between
Banach spaces. Christensen, Phelps and Mankiewicz studied the same
problem independently and used Haar null, Gaussian null and cube
null sets, respectively.
More information about the history is available in the monograph \cite{BL}.
Cs\"ornyei~\cite{C} proved that
(Borel) Aronszajn null, Gaussian null, and cube null sets coincide. It is
well known that Haar null sets form a strictly larger family than
Aronszajn null sets (see \cite{BL}).

One of the questions in differentiability theory is
to understand the structure of the sets of points of
G\^ateaux nondifferentiability of Lipschitz mappings defined on
separable Banach spaces.
The strongest result in this context is due to Preiss
and Zaj\'\i\v{c}ek~\cite{PZ}.
Let $\cA$ denote the family of Borel Aronszajn-null sets (to be defined in the sequel).
Preiss and Zaj\'\i\v{c}ek introduced a Borel $\sigma$-ideal $\tilde\cA$ such that
$\tilde\cA\sbst\cA$. It follows from a recent result of Preiss that $\cA=\tilde\cA$ in $\R^2$,
and it is unknown for $2<\dim X<\infty$. In infinite dimensions, the inclusion $\tilde\cA\sbst\cA$
is strict.
It is also unknown whether, according to the definitions of \cite{PZ},
$\tilde\cA=\tilde\cC$ and $\tilde\cC\sbst\cA^*$.

Understanding the structure of the
sets of points of non-differentiabi\-lity could possibly
also be helpful in answering the longstanding open
problem whether two separable Lipschitz isomorphic spaces
are actually linearly isomorphic. This is known for
some special Banach spaces, but is open for example for $\ell_1$
and $L_1$.

We introduce a game-theoretic approach to Aronszajn null sets.
One idea behind this approach is that the Aronszajn null sets are
defined as sets for which there exists a certain decomposition for
each complete sequence of directions, whereas in the game setting,
such a decomposition is being constructed while we are only given
one direction at a time.

Games have often been used in Banach spaces (see, e.g., the survey paper
\cite{surv}), and it would be interesting to see whether this new perspective can
yield interesting results which do not involve the new notions.

\section{The Aronszajn-null game}

Let $X$ be a separable Banach space (over $\R$).
The following definitions are classical:
\be
\item For a nonzero $x\in X$, $\cA(x)$ denotes the
collection of all Borel sets $A\sbst X$ such that
for each $y\in X$, $A\cap(\R x+y)$ has Lebesgue (one dimensional)
measure zero.
\item\label{Anull}
A Borel set $A\sbst X$ is \emph{Aronszajn-null} if for each
dense sequence $\seq{x_n}\sbst X$, there exist elements $A_n\in\cA(x_n)$,
$n\in\N$, such that $A\sbst\U_nA_n$.
\item $\cA$ denotes the collection of Aronszajn-null sets.
\ee
$\cA$ is a Borel $\sigma$-ideal.

\begin{rem}
Replacing ``dense'' by ``complete'' in item \ref{Anull}
of the above definition of Aronszajn-null sets
(i.e., requiring just that the \emph{linear span} of $\{x_n : n\in\N\}$ is dense in $X$),
one gets an equivalent definition \cite[Corollary 6.30]{BL}.
\end{rem}

The definition of Aronszajn-null sets motivates the following.
\begin{defn}
The \emph{Aronszajn-null game $\GAME{\cA}$} for a Borel set $A\sbst X$
is a game between two players,
$\I$ and $\II$, who play an inning per each natural number.
In the $n$th inning, $\I$ picks $x_n\in X$, and $\II$ responds by
picking $A_n\in\cA(x_n)$.
This is illustrated in the following figure.

\medskip

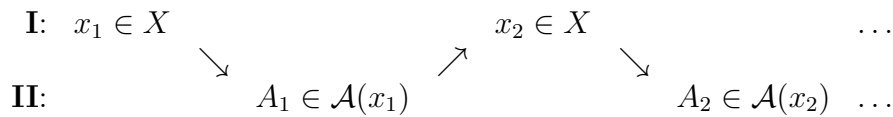
\begin{center}
\begin{tabular}{rccccccccc}
$\I$:  & $x_1\in X$ &    &      &      & $x_2\in X$ &     & & \dots\\
     &     & $\searrow$    &     & $\nearrow$    &      & $\searrow$\\
$\II$: &   &    & $A_1\in\cA(x_1)$ &   &     &    & $A_2\in\cA(x_2)$ & \dots\\
\end{tabular}
\end{center}

\medskip

\noindent $\I$ is required to play such that $\seq{x_n}$ is dense in $X$.
$\II$ \emph{wins} the game if  $A\sbst\U_nA_n$; otherwise $\I$ wins.
\end{defn}

For a game $\mathsf{G}$, the notation $\wins{\I}{\mathsf{G}}$ is a shorthand for
``$\I$ has a winning strategy in the game $\mathsf{G}$'', and
$\nwins{\I}{\mathsf{G}}$ stands for ``$\I$ does not have a winning strategy in the game $\mathsf{G}$''.
Define $\wins{\II}{\mathsf{G}}$ and $\nwins{\II}{\mathsf{G}}$ similarly.
The following is easy to see.

\begin{lem}
If $\nwins{\I}{\GAME{\cA}}$ for $A$, then $A$ is Aronszajn-null.
\hfill\qed
\end{lem}

The converse is open.

\begin{conj}\label{conj1}
If $A$ is Aronszajn-null, then $\nwins{\I}{\GAME{\cA}}$ for $A$.
\end{conj}

\begin{lem}\label{sigmaadd}
The property $\wins{\II}{\GAME{\cA}}$ is preserved under taking Borel subsets and countable unions,
i.e., it defines a Borel $\sigma$-ideal.
\end{lem}
\begin{proof}
It is obvious that $\wins{\II}{\GAME{\cA}}$ is preserved under taking Borel subsets.
To see the remaining assertion, assume that $B_1,B_2,\dots$ all satisfy $\wins{\II}{\GAME{\cA}}$,
and for each $k$ let $F_k$ be a winning strategy for $\II$ in the game $\GAME{\cA}$ played on $B_k$.
Define a strategy $F$ for $\II$ in the game $\GAME{\cA}$ played on $\Union_k B_k$
as follows.
Assume that $\I$ played $x_1\in X$ in the first inning.
For each $k$ let $A_{k,1}=F_k(x_1)$, and set $A_1=\Union_k A_{k,1}\in\cA(x_1)$.
$\II$ plays $A_1$.
In the $n$th inning we have $(x_1,A_1,x_2,A_2,\dots,x_n)$ given, where $x_n$ is
the $n$th move of $\I$.
For each $k$ let $A_{k,n}=F_k(x_1,A_{k,1},x_2,A_{k,2},\dots,x_n)$, and set $A_n=\Union_n A_{k,n}\in\cA(x_n)$.
$\II$ plays $A_n$.

Consider the play $(x_1,A_1,x_2,A_2,\dots)$.
For each $k$, $(x_1,A_{k,1},x_2,A_{k,2},\allowbreak\dots)$ is a play according to the strategy $F_k$,
and therefore $B_k\sbst\Union_n A_{k,n}$.
Consequently,
$$B=\Union_{k\in\N}B_k \sbst \Union_{k\in\N}\Union_{n\in\N} A_{k,n} = \Union_{n\in\N}\Union_{k\in\N} A_{k,n} =
\Union_{n\in\N} A_n,$$
thus $\II$ won the play.
\end{proof}

A Borel set $A\sbst X$ is \emph{directionally-porous} if there exist $\lambda>0$ and a nonzero $v\in X$
such that for each $a\in A$ and each positive $\epsilon$, there is $x\in \R v+a$
such that $\|x-a\|<\epsilon$ and $A\cap B(x,\lambda\|x-a\|)=\emptyset$.
If $A$ is directionally-porous, then so is $\cl{A}$.
$A$ is \emph{$\sigma$-directionally-porous} if it is a countable union of directionally-porous sets.

\begin{prop}\label{sdp}
For each $\sigma$-directionally-porous set, $\wins{\II}{\GAME{\cA}}$.
\end{prop}
\begin{proof}
By Lemma \ref{sigmaadd}, it suffices to consider the case where
$A\sbst X$ is directionally-porous. Let $\lambda>0$ and $v\in X$ be witnesses for that.
In this case, the function
$$F(x_1,A_1,x_2,A_2,\dots,x_n) =
\begin{cases}
A & \|x_n-v\|<\lambda/2\\
\emptyset & \mbox{otherwise}
\end{cases}$$
is a winning strategy for $\II$ in the game $\GAME{\cA}$.
\end{proof}

For a nonzero $x\in X$ and a positive $\epsilon$, let $\cA(x,\epsilon)$ denote the collection
of all Borel sets $A\sbst X$ such that for each $v\in X$ with $\|v-x\|<\epsilon$, $A\in\cA(v)$.
$\cC^*$ is the collection of all countable unions of sets $A_n$ such that each $A_n\in\cA(x_n,\epsilon_n)$
for some $x_n,\epsilon_n$. $\cC^*$ is a Borel $\sigma$-ideal.

The proof of Proposition \ref{sdp} actually establishes the following.
\begin{prop}
For each $A\in\cC^*$, $\wins{\II}{\GAME{\cA}}$.\hfill\qed
\end{prop}

The following diagram summarizes our knowledge thus far:
$$\sigma\mbox{-directionally-porous}
\Longrightarrow \cC^*
\Longrightarrow \wins{\II}{\GAME{\cA}}
\Longrightarrow \nwins{\I}{\GAME{\cA}}
\Longrightarrow \cA.$$
The open problems concerning this diagram are whether any of the last three arrows can be reversed (i.e., turned
into an equivalence) and therefore produce a characterization.
The first arrow is not reversible \cite{PZ}.

We conjecture that $\wins{\II}{\GAME{\cA}}$ is strictly stronger than $\cA$.
For brevity, we introduce the following.
\begin{defn}
For $Y\sbst X$, $A\in\cA(\wedge\ Y)$ means:
For each $y\in Y$, $A\in\cA(y)$. In other words,
$$\cA(\wedge\ Y)=\bigcap_{y\in Y}\cA(y).$$
\end{defn}
Thus, $\cA(x,\epsilon)=\cA(\wedge\ B(x,\epsilon))$.
Using this notation, we can see that the property $\wins{\II}{\GAME{\cA}}$
implies something quite close to $\cA^*$, see Corollary \ref{almostA*}.

Recall that for a topological space $X$, a pseudo-base is a family $\cU$ of open
subsets of $X$, such that each open subset of $X$ contains some element of $\cU$
as a subset. Clearly, every base is a pseudo-base.

\begin{thm}
Assume that $\wins{\II}{\GAME{\cA}}$ holds for $A$.
Then: For each countable dense $D\sbst X$ and each pseudo-base $\seq{U_n}$ for the
topology of $X$, there exist elements
$$A_n\in\cA(\wedge\ D\cap U_n),$$
$n\in\N$, such that $A\sbst\Union_n A_n$.
\end{thm}
\begin{proof}
Assume that $D\sbst X$ is countable and dense,
and $\seq{U_n}$ is a pseudo-base for the topology of $X$.
For each $n$, fix  an enumeration $\{x_{n,m} : m\in\N\}$ of $D\cap U_n$.

Let $F$ be a winning strategy for $\II$ in the game $\GAME{\cA}$.
To each finite sequence $\eta$ of natural numbers we associate
a Borel set $A_\eta$ and an element $y_\eta\in D\cap U_n$ where $n$ is the length of the sequence.
This is done by induction on $n$.

\paragraph{$n=1$:} For each $k$, set $A_k=F(x_{1,k})$.

\paragraph{$n=m+1$:} For each $\eta\in\N^m$ and each $k$, define
$$A_{\eta\cat k}=F(x_{1,\eta_1},A_{\eta|1},x_{2,\eta_2},A_{\eta|2},\dots,x_{m,\eta_m},A_\eta,x_{m+1,k}),$$
where for each $i$, $\eta_i$ is the $i$th element of $\eta$ and $\eta|i$ is the sequence $(\eta_1,\dots,\eta_i)$.

Next, for each $\eta$, define $B_\eta = \bigcap_k A_{\eta\cat k}$.
Assume that $A\not\sbst\Union_\eta B_\eta$, and let $a\in A\sm\Union_\eta B_\eta$.
Choose inductively $k_1$ such that $a\nin A_{k_1}$, $k_2$ such that
$a\nin A_{(k_1,k_2)}$, etc. Then the play $(x_{1,k_1},A_{k_1},x_{2,k_2},A_{(k_1,k_2)},\dots)$
is according to the strategy $F$ and lost by $\II$, a contradiction.
Consequently, $A\sbst\Union_\eta B_\eta$.

For each $m$ and each $\eta\in\N^m$, $B_\eta=\bigcap_k A_{\eta\cat k}\in\cA(\wedge\ D\cap U_m)$.
Consequently, $C_m = \Union_{\eta\in\N^m}B_\eta\in\cA(\wedge\ D\cap U_m)$ too, and
$A\sbst\Union_m C_m$ as required.
\end{proof}

\begin{cor}\label{almostA*}
Assume that $\wins{\II}{\GAME{\cA}}$ holds for $A$.
Then: For each countable dense $D\sbst X$, there
exist elements
$$A_n\in\cA(\wedge\ D\cap B(x_n,\epsilon_n)),$$
where each $x_n\in X$ and each $\epsilon_n>0$,
such that $A\sbst\Union_n A_n$.\hfill\qed
\end{cor}

\begin{prob}
Is the property in Corollary \ref{almostA*} equivalent to $\wins{\II}{\GAME{\cA}}$,
or does it at least imply $\nwins{\I}{\GAME{\cA}}$?
\end{prob}

\section{Selection hypotheses}

\begin{defn}
$\bbA$ is the collection of Borel sets $A\sbst X$ such that:
For each sequence $\seq{D_n}$ of dense subsets of $X$, there exist
elements $x_n\in D_n$ and $A_n\in\cA(x_n)$, $n\in\N$, such that $A\sbst\Union_n A_n$.
$\GAME{\bbA}$ is the corresponding game, played as follows:

\medskip

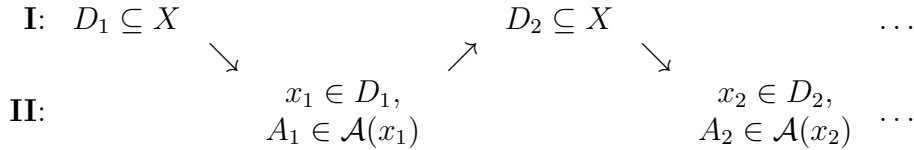
\begin{center}
\begin{tabular}{rccccccccc}
$\I$:  & $D_1\sbst X$ &    &      &      & $D_2\sbst X$ &     & & \dots\\
     &     & $\searrow$    &     & $\nearrow$    &      & $\searrow$\\
$\II$: &   &    & $\begin{matrix}x_1\in D_1,\\A_1\in\cA(x_1)\end{matrix}$ &   &     &    &
$\begin{matrix}x_2\in D_2,\\A_2\in\cA(x_2)\end{matrix}$ & \dots\\
\end{tabular}
\end{center}

\medskip
\noindent where each $D_n$ is dense in $X$, and
$\II$ wins the game if  $A\sbst\U_nA_n$; otherwise $\I$ wins.
\end{defn}

The appealing property in the game $\GAME{\bbA}$ is that, unlike the case in the
game $\GAME{\cA}$, there is no commitment of $\I$ which has to be verified
``at the end'' of the play.

\begin{prop}\label{A=A}
$\bbA=\cA$.
\end{prop}
\begin{proof}
$(\sbst)$ Assume that $A\in\bbA$, and let $D=\seq{x_n}$ be dense in $X$.
For each $n$, take $D_n=D$ and apply $\bbA$.
Then there are $y_n\in D$ and $A_n\in\cA(y_n)$, $n\in\N$, such that $A\sbst\Union_nA_n$.
As each $\cA(x_n)$ is $\sigma$-additive, we may assume that no $x_n$ appears more than
once in the sequence $\seq{y_n}$. Thus, $A\in\cA$.

$(\spst)$ Assume that $A\in\cA$, and let $\seq{D_n}$ be a sequence of dense subsets of $X$.
For each $n$ choose $x_n\in D_n$ such that $D=\seq{x_n}$ is dense in $X$ (to do that, fix
a countable base $\seq{U_n}$ for the topology of $X$, and for each $n$ pick $x_n\in U_n\cap D_n$).
By $\cA$, there exist sets $A_n\in\cA(x_n)$ such that $A\sbst\Union_nA_n$. This shows that $A\in\bbA$.
\end{proof}

A simple modification of the last proof gives the following.
\begin{thm}\label{GameA=A}
$\wins{\I}{\GAME{\bbA}}$ if, and only if, $\wins{\I}{\GAME{\cA}}$.
\end{thm}
\begin{proof}
$(\Impl)$ Let $F$ be a winning strategy for $\I$ in the game $\wins{\I}{\GAME{\bbA}}$ on $A$.
Define a strategy for $\I$ in the game $\wins{\I}{\GAME{\cA}}$ as follows.
Fix a countable base $\seq{U_n}$ for the topology of $X$.
In the first inning, $\I$ plays any $x_1\in U_1\cap D_1$ where $D_1$ is $\I$'s first move
according to the strategy $F$.
Assume that the first $n$ moves where $(x_1,A_1,\dots,x_{n-1},A_{n-1})$.
Let $D_n = F(D_1,(x_1,A_1),\dots,D_{n-1},(x_{n-1},A_{n-1}))$.
Then $\I$ plays any $x_n\in U_n\cap D_n$.
For each play $(x_1,A_1,x_2,A_2,\dots)$ according to this strategy, $\seq{x_n}$ is dense in $X$, and
since $(D_1,(x_1,A_1),D_2,(x_2,A_2),\dots)$ is a play in the game $\GAME{\bbA}$ according to the strategy $F$,
$A\not\sbst\Union_n A_n$.

$(\Leftarrow)$ Let $F$ be a winning strategy for $\I$ in the game $\wins{\I}{\GAME{\cA}}$ on $A$.
Define a strategy for $\I$ in the game $\wins{\I}{\GAME{\bbA}}$ as follows.
$\I$'s first move is $D_1$, the set of all points $x$ which are possible moves of $\I$ at some inning
according to its strategy $F$. Obviously, $D_1$ is dense.
In the $n$th inning, we are given $(D_1,(x_1,A_1),\dots,D_{n-1},(x_{n-1},A_{n-1}))$, such
that there is a sequence of moves $(y_1,B_1,y_2,B_2,\dots,y_{k_n})$ according to the strategy
$F$, with $y_{k_n}=x_{n-1}$.
Then $\I$ plays $D_n$, the set of all points $x$ which are possible moves of $\I$ at some future inning,
in a play according to the strategy $F$ whose first moves are $(y_1,B_1,y_2,B_2,\dots,y_{k_n}=x_{n-1},A_{n-1})$.
$(y_1,B_1,y_2,B_2,\dots)$ is a play according to the strategy $F$, and therefore
$A\not\sbst\Union_n B_n\spst\Union_n A_n$, so that $A\not\sbst\Union_n A_n$.
\end{proof}

The following is immediate.

\begin{lem}
If $\wins{\II}{\GAME{\cA}}$, then $\wins{\II}{\GAME{\bbA}}$.\hfill\qed
\end{lem}

\begin{prob}
Is it true that $\wins{\II}{\GAME{\bbA}}$ if, and only if, $\wins{\II}{\GAME{\cA}}$?
\end{prob}

By Theorem \ref{GameA=A} and Proposition \ref{A=A}, Conjecture \ref{conj1} can be reformulated
as follows.

\begin{conj}
For a Borel set $A\sbst X$: $\nwins{\I}{\GAME{\bbA}}$ if, and only if, $A\in\bbA$.
\end{conj}

\subsection*{Acknowledgments}
The second author was partially supported by the Koshland Center for Basic Research.
We thank Ori Gurel-Gurevich for presenting our results at the Weizmann Institute's
seminar on Geometric Functional Analysis and Probability, and for making useful comments.

\ed